\documentclass[letterpaper, 10 pt, conference]{ieeeconf}  

\IEEEoverridecommandlockouts                              

\usepackage{cite}
\usepackage{amsmath,amssymb,amsfonts}
\usepackage{algorithmic}
\usepackage[ruled]{algorithm2e} 
\usepackage{graphicx}
\usepackage{textcomp}
\usepackage{xcolor}
\usepackage{balance}

\newtheorem{theorem}{Theorem}
\newtheorem{definition}{Definition}
\newtheorem{remark}{Remark}
\newtheorem{assumption}{Assumption}

\newtheorem{lemma}{Lemma}
\newtheorem{proposition}{Proposition}

\def\BibTeX{{\rm B\kern-.05em{\sc i\kern-.025em b}\kern-.08em
    T\kern-.1667em\lower.7ex\hbox{E}\kern-.125emX}}

\title{\LARGE \bf
Incremental Policy Iteration for Unknown Nonlinear Systems with Stability and Performance Guarantees
}

\author{Qingkai Meng, Fenglan Wang, and Lin Zhao
\thanks{All authors are with Department of Electrical and Computer Engineering, National University of Singapore, 117583, Singapore.
        {\tt\small qk.meng@nus.edu.sg; wfenglan@nus.edu.sg; elezhli@nus.edu.sg.} Correspondence to: Lin Zhao.}
\thanks{This work was supported by the Singapore Ministry of Education Tier~1 Academic Research Fund (A-8001174-00-00) and Tier~2 Academic Research Funds (T2EP20123-0037).}
}

\begin{document}


\maketitle

\begin{abstract}
This paper proposes a general incremental policy iteration adaptive dynamic programming (ADP) algorithm for model-free robust optimal control of unknown nonlinear systems. 
The approach integrates recursive least squares estimation with linear ADP principles, which greatly simplifies the implementation while preserving adaptive learning capabilities. In particular, we develop a sufficient condition for selecting a discount factor such that it allows learning the optimal policy starting with an initial policy that is not necessarily stabilizing. 
Moreover, we characterize the robust stability of the closed-loop system and the near-optimality of iterative policies. 
Finally, we perform numerical simulations to demonstrate the effectiveness of the proposed method.
\end{abstract}


\section{Introduction}
Adaptive Dynamic Programming (ADP) is a powerful method for solving the Hamilton-Jacobi-Bellman (HJB) equation in optimal control problems for uncertain and nonlinear systems~\cite{Bertsekas2012dynamic,Wang2024recent}. It is particularly effective in mitigating the curse of dimensionality by approximating the value function and control policy using function approximators, such as neural networks and polynomial functions.
ADP can be categorized into two primary algorithmic frameworks: policy iteration (PI) and value iteration (VI). Both approaches typically rely on high-fidelity simulation models in offline training. However, obtaining such models can be challenging in practical applications involving highly nonlinear and uncertain physical systems~\cite{Lewis2013reinforcement}.

For a linear time-invariant system, the system matrices are inherently embedded in the temporal evolution of input-output data, as they satisfy an overdetermined equation constraint. This allows the system model to be fully represented by data, enabling linear ADP methods to achieve optimal control without requiring an explicit system model~\cite{Kiumarsi2015optimal,Gao2016adaptive}.  
A fundamental question arises:  \emph{Can nonlinear systems, like their linear counterparts, leverage data-driven strategies to achieve optimal control}, thereby eliminating the need for explicit model identification or neural network model training? 
Addressing this challenge requires extending linear ADP methods to nonlinear system control.
One promising approach is incremental control~\cite{Sieberling2010robust}, which iteratively refines control strategies using locally linearized approximations.   
Building on this idea, recent research has explored a hybrid framework that integrates nonlinear incremental control techniques with linear VI-ADP~\cite{Zhou2017nonlinear}.
While this method has demonstrated effectiveness in numerical simulations~\cite{Meng2022fault}, its stability guarantees remain underdeveloped, highlighting the need for further theoretical analysis. The importance of stability-guaranteed design in uncertain systems has been widely recognized in adaptive control research. A representative work is~\cite{Li2024Automatica}, which demonstrates closed-loop stability and asymptotic tracking for uncertain nonlinear systems.

Although stability guarantees for standard ADP methods are well-established in both linear and nonlinear settings~\cite{Vrabie2009adaptive,Bian2016value,Liu2013Policy,Wei2015value}, they can not directly apply to ADP algorithms with incremental models. 
The key challenge lies in the model approximation errors introduced by using an incremental linear model to represent nonlinear dynamics, which complicates convergence and stability analysis.
Furthermore, since VI does not explicitly refine the control policy at each iteration, it lacks a direct mechanism to ensure improvement in stability during training, making VI less suitable for online learning~\cite{Vrabie2009adaptive}.  
Although PI can provide online learning solutions with stability guarantees, it requires an initially stabilizing policy, which is computationally expensive to obtain, especially when the system model is unknown~\cite{Granzotto2024robust}.
Given these challenges, developing incremental PI-ADP algorithms without initially stabilizing policy requirement is of great significance for achieving online learning control with theoretical guarantees. 

Motivated by the above discussions, this paper proposes an Incremental Policy Iteration (IPI) framework, which reformulates the integration of the incremental control technique and PI while eliminating the need for an initially stabilizing policy. 
Theoretical guarantees for IPI, including near-optimality and robust stability, are rigorously established. 
The main contributions of this work are threefold:
\begin{itemize}
    \item[i)] A general IPI algorithm is proposed, where a first-order Taylor series approximation model of the system dynamics is identified using recursive least squares (RLS) methods~\cite{Isermann2011introduction}.
    This facilitates next-step state computation during  offline  training while enabling online policy optimization with limited data, providing a model-free controller design adaptable to dynamic system variations.
    \item[ii)] To handle errors from linearized approximations, system identification, and value-function surrogates, we design iteration rules that synchronize model updates, value approximations, and policy improvements. These rules bound the accumulated errors, ensuring that the generated policies remain provably near-optimal for nonlinear control using linear ADP methods.
    \item[iii)] Using a general indicator function to define attractor-related stability, we show via Lyapunov analysis that the proposed IPI algorithm converges to a stabilizing policy, given a sufficient number of iterations. Furthermore, we establish an explicit relationship between closed-loop stability and the discount factor, relaxing the requirement for an initially stabilizing policy.
\end{itemize}


\section{Preliminaries}\label{sec:preliminaries}
\subsection{Notations}
Denote the set of real numbers by $\mathbb{R}$, the set of integers by $\mathbb{Z}$, and the set of $n$ dimensional real number vectors by $\mathbb{R}^n$. 
Denote a subset of $\mathcal{A}$ satisfying $(\cdot)$ by $\mathcal{A}_{(\cdot)}$.
A function \(\alpha: [0, a) \to [0, \infty)\) is of \emph{class \(\mathcal{K}\)} if \(\alpha(0) = 0\) and \(\alpha(r)\) is continuous and strictly increasing; it is of \emph{class \(\mathcal{K}_{\infty}\)} if additionally \(a = \infty\) and \(\alpha(r) \to \infty\) as \(r \to \infty\); a function \(\beta: [0, \infty) \times [0, \infty) \to [0, \infty)\) is of \emph{class \(\mathcal{KL}\)} if for each fixed \(t \geq 0\), \(\beta(r, t)\) is of class \(\mathcal{K}\), and for each fixed \(r \geq 0\), \(\beta(r, t)\) is continuous, strictly decreasing, and \(\beta(r, t) \to 0\) as \(t \to \infty\).
The Euclidean norm of a vector $x\in\mathbb{R}^{n_x}$ with $n_x\in\mathbb{Z}_{>0}$ is denoted by $\|x\|$ and the distance of $x\in\mathbb{R}^{n_x}$ to a nonempty closed set $\mathcal{A}\subset\mathbb{R}^{n_x}$ is denoted by $\|x\|_{\mathcal{A}}:\inf\{\|x-y\|:y\in\mathcal{A}\}$.
Given $\mathcal{A}$, the map $\sigma(x):\mathbb{R}^{n_x}\rightarrow \mathbb{R}_{\geq 0}$ is a \emph{proper indicator function} of set $\mathcal{A}$ whenever $\sigma$ is continuous and there exist $\underline{\sigma},\bar{\sigma}\in\mathcal{K}_{\infty}$ such that $\underline{\sigma}(\|x\|_{\mathcal{A}})\leq \sigma(x)\leq \bar{\sigma}(\|x\|_{\mathcal{A}})$.

\subsection{Model and cost function}
Consider a nonlinear system in the form of
\begin{equation}\label{eq:nonlinear system}
    x_{k+1}=f(x_k,u_k),\ \forall k\in\mathbb{Z}_{\geq 0},
\end{equation}
where $x_k:=x(t_k)\in\Omega\subset\mathbb{R}^{n_x}$ is the state on the compact set $\Omega$, $u_k:=u(t_k)\in\mathcal{U}(x_k)\subseteq\mathbb{R}^{n_u}$ 
is the control input at time instant $t_k$ with $t_{k+1}=t_k+\Delta t$, $\Delta t>0$ is a fixed sampling time interval,
$\mathcal{U}(x_k)$ is a non-empty compact set of admissible inputs at state $x_k$, $k\in\mathbb{Z}_{\geq 0}$, and $n_x,n_u\in\mathbb{Z}_{>0}$ are the dimensions of state and control input, respectively.
The vector field $f(\cdot,\cdot):\mathbb{R}^{n_x}\times\mathbb{R}^{n_u}\rightarrow \mathbb{R}^{n_x}$ is unknown but assumed to be Jacobian-Lipschitz
on $\Omega\times\mathcal{U}$

The solution to \eqref{eq:nonlinear system} is denoted 
 by $\phi(k,x,{\mathbf{u}\mid_k})$ at time $t_k$ with the initial state $x$ and an admissible truncated control sequence ${\mathbf{u}\mid_k}:=\{u(0),\cdots,u(k-1)\}$. 
 We use the convention $\phi(0,x,{\mathbf{u}\mid_0})=x$. 
 We wish to find an infinite-length sequence of admissible inputs $\mathbf{u}$ by using the available data that minimizes the infinite horizon cost
 \begin{equation}\label{eq:cost function}
 \begin{aligned}
     J_{\gamma}(x,\mathbf{u}):=&\sum_{k=0}^{\infty} \gamma^k\ell\big(\phi(k,x,{\mathbf{u}\mid_k}),u_k\big),
\end{aligned}
 \end{equation}
 where $\gamma\in(0,1)$ is a cost discount factor, and $\ell:\mathbb{R}^{n_x}\times\mathbb{R}^{n_u}\rightarrow \mathbb{R}_{\geq 0}$ is a non-negative stage cost such that $|\ell(x,u)-\ell(y,u)|\leq L_{\ell}\|x-y\|$ with a known constant $L_{\ell}>0$, $x,y\in\mathbb{R}^{n_x}$.
 For any $x\in\Omega$, the optimal value function associated with the minimization of \eqref{eq:cost function} is denoted by
 \begin{equation}\label{eq:optimal value function}
     V_\gamma^{\star}(x):=\min_{\mathbf{u}}J_{\gamma}(x,\mathbf{u})<+\infty.
 \end{equation}
 As a result, the Bellman equation becomes
 \begin{equation*}
     V_\gamma^{\star}(x)=\min_{u\in\mathcal{U}}\{\ell(x,u)+\gamma V_\gamma^{\star}(f(x,u))\}, \forall x\in\Omega.
 \end{equation*}
 The optimal inputs for any state $x\in\Omega$ constitute a non-empty set as
 \begin{equation}\label{eq:optimal control set}
     H^{\star}_\gamma(x):=\arg\min_{u\in\mathcal{U}}\{\ell(x,u)+\gamma V_\gamma^{\star}(f(x,u))\}.
 \end{equation}
For the nonlinear system~\eqref{eq:nonlinear system} with a general cost function~\eqref{eq:cost function}, computing $H^{\star}_\gamma$ in \eqref{eq:optimal control set} is extremely challenging, especially when the system dynamics are unknown.
Therefore, it is necessary to utilize dynamic programming iterations to obtain the feedback law, ensuring that its cost asymptotically converges to the optimal value.

To achieve this, the following necessary assumptions are given~\cite{Bertsekas2012dynamic}.
The existence of the optimal and stabilizing control sequence, as defined below, is a prerequisite for optimization.
\begin{assumption}\label{as:existense of optimal control}
   For any $x\in \mathbb{R}^{n_x}$ and any $\gamma\in (0,1)$, there exists an optimal sequence of admissible inputs $\mathbf{u}^{\star}(x)$ such that $V^{\star}_{\gamma}(x)=J_{\gamma}(x,\mathbf{u}^{\star}(x))<\infty$ and for any infinite-length sequence of admissble inputs $\mathbf{u}$, $V_{\gamma}^{\star}(x)\leq J_{\gamma}(x,\mathbf{u})$.  \hfill$\square$
\end{assumption}


\begin{assumption}\label{as:stabilizing optimal policy}
    There exists $\bar{\alpha}_{V^{\star}}\in\mathcal{K}_{\infty}$ and $\gamma_0\in(0,1]$ such that for any $\gamma\in(0,\gamma_0)$ and any $x\in\Omega$, $V^{\star}_{\gamma}(x)\leq \bar{\alpha}_{V^{\star}}(\sigma(x))$, where $V^{\star}_{\gamma}$ is given in \eqref{eq:optimal value function}.
      \hfill$\square$
\end{assumption}

\section{Problem Formulation}\label{sec:problem formulation}

\subsection{Incremental policy iteration}
This subsection introduces the incremental policy iteration, as shown in Algorithm~\ref{al PI-IADP}, to iteratively obtain feedback laws.

Specifically, taking the Taylor expansion of \eqref{eq:nonlinear system} at state $x_k$, the following incremental model is obtained
\begin{equation}\label{eq:incremental model}
    \Delta x_{k+1}=A_{k-1}\Delta x_k+B_{k-1}\Delta u_k+O(\Delta x_k^2, \Delta u_k^2),
\end{equation}
where the incremental state and control at time $k $ are defined as $\Delta x_{k}:=x_k-x_{k-1}$ and $\Delta u_k=u_k-u_{k-1}$, matrices $A_{k-1}:=\frac{\partial f}{\partial x}\mid_{x=x_{k-1}}$ and $B_{k-1}:=\frac{\partial f}{\partial u}\mid_{u=u_{k-1}}$ are partial derivatives of the dynamics with respect to the state and control at time $t_{k-1}$, and $O(\cdot)$ denotes the high-order remainder.
Since $f$ is Jacobian-Lipschitz, the $O(\cdot)$ is bounded on $\Omega\times\mathcal{U}$.
The nonlinear system \eqref{eq:nonlinear system} can be represented as this time-varying incremental model, in which $A_{k-1},B_{k-1}$ are expected to be identified by using RLS methods~\cite{Isermann2011introduction}.

\begin{algorithm}[tb]\label{al PI-IADP}
\caption{Incremental Policy Iteration}
\KwIn{State $x_k,x_{k-1}$, initial policy $h_{\gamma}^0\in\mathcal{U}$, initial system matrices $\hat{\Theta}_0=[\hat{A}_0,\hat{B}_0]^{\top}$, initial covariance matrix $\Lambda_0$, RLS discounted factor $\kappa$, stage cost $\ell(\cdot,\cdot)$, initial approximator $W^0_{\gamma}$}
\KwOut{Policy $u_{\gamma}^{\infty}$, cost $V_{\gamma}^{\infty}$}
\begin{algorithmic}[1]
\raggedright
\STATE {\textbf{RLS Identification:}\\
1.1: $\Delta\hat{x}_{k+1}^{\top}=X_k^{\top}\hat{\Theta}_{k-1}$, $X_k:=[\Delta x_k; \Delta u_k]\makebox[43pt][r]{(LS.1)}$\\
1.2: $\varepsilon_k=\Delta x_{k+1}^{\top}-\Delta\hat{x}_{k+1}^{\top}\makebox[114pt][r]{(LS.2)}$\\
1.3: $\hat{\Theta}_k=\hat{\Theta}_{k-1}+\frac{\Lambda_{k-1}X_k}{\kappa+X^{\top}_k\Lambda_{k-1}X_k}\varepsilon_k\makebox[81pt][r]{(LS.3)}$\\
1.4: $\Lambda_k=\frac{1}{\kappa}\big[\Lambda_{k-1}-\frac{\Lambda_{k-1}X_kX_k^{\top}\Lambda_{k-1}}{\kappa+X_k^{\top}\Lambda_{k-1}X_k}\big]\makebox[70pt][r]{(LS.4)}$\\
1.5: Return $\hat{A}_k,\hat{B}_k$}
\STATE {\textbf{Policy Iteration:}\\
~~2.1: \textbf{Initial evaluation step:}\\ 
~~~~~~$\hat{x}_{k+1}=x_k+\hat{A}_{k-1}\Delta x_k+\hat{B}_{k-1}\Delta h^0_{\gamma,k}\makebox[45pt][r]{(PI.1)}$\\
~~~~~~~~$\hat{V}_{\gamma}^0(x_k):=\ell(x_k,h^{0}_{\gamma})+\gamma W^0_{\gamma}(\hat{x}_{k+1})\makebox[51pt][r]{(PI.2)}$\\
~\textbf{for} $i\in\mathbb{Z}_{\geq 0}$ \textbf{do}\\
~~2.2: \textbf{Policy improvement step:}\\ ~~~~~~$\hat{x}_{k+1}=x_{k}+\hat{A}_{k-1}\Delta x_k+\hat{B}_{k-1}\Delta h^i_{\gamma,k}\makebox[45pt][r]{(PI.3)}$\\
~~~~~~~$\hat{V}_{\gamma}^{i}(x_k):=\ell(x_k,h^{i}_{\gamma})+\gamma W^{i}_{\gamma}(\hat{x}_{k+1})\makebox[56pt][r]{(PI.4)}$
\\
~~~~~~~$\Delta H_\gamma^{i+1}(x_{k}):=\mathop{\arg\min}\limits_{\Delta u_k\in\Delta\mathcal{U}(x_k)}\hat{V}^{i}_{\gamma}(x_k)
\makebox[56pt][r]{(PI.5)}$\\
~~\textbf{Select} $\Delta h_{\gamma}^{i+1}\in \Delta H_{\gamma}^{i+1}$, $h_{\gamma}^{i+1}=h_{\gamma}^{i}\!+\!\Delta h_{\gamma}^{i+1}\!\in\! H_{\gamma}^{i+1}$\\
~~2.3: \textbf{Policy evaluation step:}\\ 
~~~$W_{\gamma}^{i+1}(x_k)=\ell(x_k,h_{\gamma}^{i+1})+\gamma W^{i+1}_{\gamma}(\hat{x}_{k+1})\makebox[42pt][r]{(PI.6)}
$\\
~\textbf{end for}\\
\textbf{Return} $u_{\gamma}^{\infty}\in H_{\gamma}^{\infty}$ and $V_{\gamma}^{\infty}$}
\end{algorithmic}
\end{algorithm}

The concrete RLS identification principle is given below. 
With the following augmented system state 
\begin{equation*}
    X_k:=\begin{bmatrix}
        \Delta x_k\\
        \Delta u_k
    \end{bmatrix},
\end{equation*}
the augmented system matrices 
\begin{equation*}
    \hat{\Theta}_{k-1}:=\big[\hat{A}_{k-1}\ \hat{B}_{k-1}\big]^{\top}
\end{equation*}
are responsible for the one-step prediction as 
\begin{equation*}
    \Delta \hat{x}^{\top}_{k+1}=X_k^{\top}\hat{\Theta}_{k-1}.
\end{equation*}

The identification is achieved by using ($\mathrm{LS.1}$)-($\mathrm{LS.4}$) in Algorithm~\ref{al PI-IADP} with a recursive manner.
For the physical system matrices $\Theta_k, k\in\mathbb{Z}_{\geq 0}$, denote the estimation error by $\tilde{\Theta}_k:=\hat{\Theta}_k-\Theta_k$. 
Define the system incremental error as $\Delta\Theta_k:=\Theta_k-\Theta_{k-1}$, which is determined by the dynamics~\eqref{eq:nonlinear system}. 
Based on the form of incremental model~\eqref{eq:incremental model},  
it follows that there exists an upper bound $\varepsilon_{\Delta\Theta}$ such that $\|\Delta\Theta_k\|\leq \varepsilon_{\Delta\Theta}$, for all $k\in\mathbb{Z}_{>0}$.
As a result, the estimation error is bounded and satisfies the following properties.
\begin{lemma}[\cite{Isermann2011introduction}]\label{le:RLS error}
 For system~\eqref{eq:nonlinear system}, denote the changing of incremental model~\eqref{eq:incremental model} upper bound by $\varepsilon_{\Delta\Theta}>0$.
 Then, when identifying $\Theta_k$ using (\emph{$\mathrm{LS.1}$})-(\emph{$\mathrm{LS.4}$}), there exists $\beta_{\Theta}\in\mathcal{K}\mathcal{L}$ such that $\tilde{\Theta}_k\leq \beta_{\Theta}(\varepsilon_{\Delta\Theta},k)$, $k\in\mathbb{Z}_{>0}$. \hfill$\square$
\end{lemma}

With the estimated matrices $\hat{\Theta}_k$ at time $t_k$, the estimated incremental state at successor time $t_{k+1}$ is given by
\begin{equation}\label{eq:indentified incremental model}
    \Delta\hat{x}_{k+1}= \hat{A}_{k-1}\Delta x_k+\hat{B}_{k-1}\Delta u_k.
\end{equation}
A parameterized value function approximator $W^i_{\gamma}(\cdot)$ is introduced and the approximation value $\hat{V}_{\gamma}(\hat{x}_{k+1})$ at time $k+1$ is given by substituting $\hat{\phi}(1,x_k,u_k):=\hat{x}_{k+1}=x_k+\Delta\hat{x}_{k+1}$ into the approximator $W^i_{\gamma}(\cdot)$, yielding ($\mathrm{PI.4}$). 

The IPI is assumed to satisfy the following conditions.
\begin{assumption}\label{as:IPI}
  There exist a constant $0< \gamma_0\leq 1$, functions $\bar{\alpha}_{V}(\cdot,\gamma),\alpha_\Gamma(\cdot)\in\mathcal{K}_{\infty}$, a continuous function $\Gamma:\mathbb{R}^{x_n}\rightarrow \mathbb{R}_{\geq0}$,  such that $\forall \gamma\in(0,\gamma_0]$,
  \begin{subequations}
      \begin{align}
          \ell\big(x,h_{\gamma}^0(x)\big)\!+\!\gamma W_{\gamma}^0\big(\hat{\phi}(1,x,h_{\gamma}^0(x))\big)\leq \bar{\alpha}_{V}(\sigma(x),\gamma), \label{eq1 Assumption for IPI} \\
          \begin{aligned}
          \Gamma&\big(\hat{\phi}(1,x,h_{\gamma}\big)-\Gamma(x)\leq -\alpha_\Gamma(\sigma(x))+\ell\big(x,h_{\gamma}\big),\label{eq2 Assumption for IPI}
         \end{aligned}\\
         \hat{V}^i_{\gamma}(x)\geq\ell(x,h_{\gamma}^{i+1}(x))+\gamma W^i_{\gamma}(\hat{\phi}(1,x,h_{\gamma}^{i+1}(x))),\label{eq3 Assumption for IPI}
      \end{align}
  \end{subequations}
  for all $i\in\mathbb{Z}_{\geq 0}$, $x\in\Omega$.
  \hfill$\square$
\end{assumption}

\begin{remark}
 Condition~\eqref{eq1 Assumption for IPI} establishes the relationship between the stability (the distance of state $x$ from the attractor) of the original system~\eqref{eq:nonlinear system} and the initial estimated value function $\hat{V}^{0}_{\gamma}$. 
 There exists no assumption of stabilizing initial policy but only the bounded initial estimated value function by $\bar{\alpha}_{V}$ is required.
 This is reasonable and easy to implement because, for any given state \(x \in \mathbb{R}^{n_x}\) and an initial policy \(h_{\gamma}^0 \in \mathcal{U}(x)\), \(\sigma(x)\) and \(\ell(x, h_{\gamma}^0(x))\) can be computed explicitly. 
 By choosing an appropriate \(\gamma\) and an approximation operator \(W_{\gamma}^0(\cdot)\), condition \eqref{eq1 Assumption for IPI} can be satisfied.
 In particular, the smaller the \(\gamma\), the easier it is to satisfy this condition, which can be reflected in selecting a sufficiently small \(\gamma_0\).
 An appropriate form of the initial policy $h^0$ can be selected to satisfy $\ell(\hat{\phi}(k,x,h^0),h^0(\phi(k,x,h^0)))\leq Ma^k \chi(\sigma(x))$ with $M,a>0$ and $\chi\in\mathcal{K}_{\infty}$, which may be not stabilizing when $a$ is strictly bigger than 1~\cite{de2024policy}. 

 Condition~\eqref{eq2 Assumption for IPI} 
specifies the detectability property of system~\eqref{eq:nonlinear system} with incremental model approximation~\eqref{eq:incremental model}, when consider $\ell$ as an output.
This is natural as this shows the fact that by minimizing $\ell$ along the solution to~\eqref{eq:incremental model}, desirable stability properties should follow.

Condition~\eqref{eq3 Assumption for IPI} indicates that the cost generated by the next policy \(h_{\gamma}^{i+1}\), along with the discounted future cost, should not exceed the current estimated total cost.
That is, each updated policy \(h_{\gamma}^{i+1}\) is better than or at least equivalent to the previous policy \(h_{\gamma}^i\).
This is consistent with the upper bound form of the Bellman equation and ensures stable policy improvement.
   \hfill$\square$
\end{remark}
 


With the IPI, the closed-loop system with approximation optimal controller is given by
\begin{equation}\label{eq:closed loop system}
\begin{aligned}
    x_{k+1}&\in f\big(x_k,H_{\gamma}^{i+1}(x_k)\big)\\
    &=x_k+\hat{A}_{k-1}\Delta x_k+\hat{B}_{k-1}\Delta h^{i}_{\gamma,k}+\Delta_{\mathrm{IME}},
\end{aligned}
\end{equation}
where $\Delta h^{i}_{\gamma,k}\in \Delta H_\gamma^{i}(x_k)$ defined in ($\mathrm{PI.5}$) and $\Delta_{\mathrm{IME}}:=(\hat{A}_{k-1}-A_{k-1})\Delta x_k+(B_{k-1}-\hat{B}_{k-1})\Delta u_k+O(\Delta x_k^{2},\Delta (u^i_k)^2)$ is the total error of using IPI.
From Lemma~\ref{le:RLS error} and  the boundedness of $O(\cdot)$, 
it is reasonable to assume that there exists a constant $\varepsilon_{\mathrm{IME}}>0$ such that
$\|\Delta_{\mathrm{IME}}\|\leq \varepsilon_{\mathrm{IME}}$.
Therefore, the upper bound of the error between the optimal value function $\hat{V}^{\star}_{\gamma}$ of the incremental model and that of the original nonlinear system~\eqref{eq:nonlinear system} is deduced below.
\begin{proposition}
  Consider the system~\eqref{eq:nonlinear system} with the incremental approximation~\eqref{eq:closed loop system} controlled by policies generated from Algorithm~\ref{al PI-IADP}. 
  If there exists constant $\varepsilon_{\mathrm{IME}}>0$ such that $\|\Delta_{\mathrm{IME}}\|\leq\varepsilon_{\mathrm{IME}}$, then for all $k\in\mathbb{Z}_{\geq 0}$ and $\gamma\in (0,1)$,
  \begin{equation}\label{eq:bound of incremental iotimal}
      |V_{\gamma}^{\star}(x_k)-\hat{V}_{\gamma}^{\star}(x_k)|\leq \frac{\gamma L_{\ell}\varepsilon_{\mathrm{IME}}}{1-\gamma},
  \end{equation}
  where $L_{\ell}$ is the Lipscitz constant of $\ell(\cdot,\cdot)$. \hfill$\square$
\end{proposition}

\textit{Proof.}
    For any initial state $x_0\in\mathbb{R}^{n_x}$, $\gamma\in(0,1)$, $i\in\mathbb{Z}_{\geq 0}$, solution $x_{k+1}$ to \eqref{eq:closed loop system} and solution $\hat{x}_{k+1}$ to (PI.3),  
    let $u^{\star}_k\in H^{\star}_{\gamma}$ and $\hat{u}^{\star}_k\in \hat{H}^{\star}_{\gamma}$.
    By Bellman equation and the definition of $V_{\gamma}^{\star}$,
    \begin{equation*}
    \begin{aligned}
        \hat{V}^{\star}_{\gamma}(x_0)-V^{\star}_{\gamma}(x_0)=&\sum_{k=0}^{\infty}\gamma^k\ell(\hat{x}_{k},\hat{u}^{\star}_k)\!-\!\sum_{k=0}^{\infty}\gamma^k\ell({x}_{k},{u}^{\star}_k).
    \end{aligned}
    \end{equation*}
    As the following equation holds,
    \begin{equation*}
    \hat{u}^{\star}_k\in \hat{H}^{\star}_{\gamma}(\hat{x}_{k-1})=\mathop{\arg\min}\limits_{\hat{u}_k\in\mathcal{U}(x_{k-1})}\{\ell(\hat{x}_{k-1},\hat{u}_{k-1})+\hat{V}^{\star}_{\gamma}(\hat{x}_{k+1})\},
    \end{equation*}
    we have that for all $\hat{u}_{k}\in\mathcal{U}(\hat{x}_{k-1})$, $k\in\mathbb{Z}_{>0}$,
    \begin{equation*}
       \begin{aligned}
        \hat{V}^{\star}_{\gamma}(x_0)-V^{\star}_{\gamma}(x_0)\leq&\sum_{k=0}^{\infty}\gamma^k\ell(\hat{x}_{k},\hat{u}_k)\!-\!\sum_{k=0}^{\infty}\gamma^k\ell({x}_{k},{u}^{\star}_k).
    \end{aligned}
    \end{equation*}
    When selecting $\hat{u}_k=u^{\star}_k$, it follows that
    \begin{equation*}
        \begin{aligned}
        |\hat{V}^{\star}_{\gamma}(x_0)-V^{\star}_{\gamma}(x_0)|\leq&\sum_{k=0}^{\infty}\gamma^k\big|\ell(\hat{x}_{k},{u}^{\star}_k)-\ell({x}_{k},{u}^{\star}_k)\big|\\
        \leq& \sum_{k=0}^{\infty}\gamma^k L_{{\ell}}\varepsilon_{\mathrm{IME}}=\frac{\gamma L_{\ell}\varepsilon_{\mathrm{IME}}}{1-\gamma}.
        \end{aligned}
    \end{equation*}
    This completes the proof.\hfill$\blacksquare$

Denote $\rho(x_k):=\|\Delta_{\mathrm{IME}}\|$, which inspires the following perturbation set-valued closed-loop system
\begin{equation}\label{eq:closed deffirentail inclusion}
    \begin{aligned}
    x_{k+1}\in x_k+\hat{A}_{k-1}\Delta x_k+\hat{B}_{k-1}\Delta h^i_{\gamma,k}+\rho(x_k)\mathbb{B},
\end{aligned}
\end{equation}
where $\mathbb{B}$ is the unit closed ball of $\mathbb{R}^{n_x}$ centered at the origin, and $\Delta h^i_{\gamma,k}\in \Delta H^i_{\gamma}$ is defined in (PI.5).

\subsection{Study objectives}
The definitions of desired properties for IPI are given.

\begin{definition}[Near-optimality,\cite{Granzotto2024robust}]\label{def:near optimality}
  Consider system~\eqref{eq:nonlinear system} with the infinite horizon cost~\eqref{eq:cost function} and the minimization value $V_{\gamma}^{\star}$.
  A policy iteration algorithm is with \emph{near-optimality} if there exists a bound $\varepsilon_{V^{\star}}$ such that $\hat{V}^i_{\gamma}(x)-V^{\star}_{\gamma}(x)\leq \varepsilon_{V^{\star}}$, $\forall i\in\mathbb{Z}_{\geq 0}$ and $x\in\Omega$. \hfill$\square$
\end{definition}

\begin{definition}[Robust stability,\cite{Granzotto2024robust}]\label{def:robust stability}
The system~\eqref{eq:closed loop system} with the policy $h_{\gamma}^i$, $i\in\mathbb{Z}_{\geq 0}$, is robustly stable if there exists $\beta\in\mathcal{KL}$ and any given bound $\bar{\delta}\geq 0$ such that $\sigma(\phi(k,x,h_{\gamma}^{i}))\leq\max\{\beta(\sigma(x),k),\bar{\delta}\}$ for every $x\in\Omega$ and $k\in\mathbb{Z}_{\geq0}$.  
\hfill$\square$
\end{definition}


Since it is impractical to implement Algorithm 1 with infinite iterations and the approximation errors introduced by the incremental model, RLS, and the approximator cannot be ignored, our objectives are to 
\begin{itemize}
    \item[i)] verify that it 
 provides an approximate optimality guarantee, where the value error is bounded,
    \item[ii)] and establish conditions under which the IPI can produce a robustly stable control policy within a finite number of iterations.
\end{itemize}

\section{Main Results}\label{sec:main results}

\subsection{Near-optimality}
We first give the result on the improvement property of the IPI with respect to the value function.
\begin{proposition}\label{pro:improvement property}
 Under Assumption~\ref{as:IPI}, for any $x\in\Omega$, $\gamma\in(0,\gamma_0)$, $\hat{V}_{\gamma}^{i+1}(x)\leq \hat{V}_{\gamma}^{i}(x)$, $i\in\mathbb{Z}_{\geq 0}$.  \hfill$\square$
\end{proposition}

\textit{Proof.}
    Evaluating \eqref{eq3 Assumption for IPI} at $x_0\in\mathbb{R}^{n_x}$, one has
    \begin{equation}\label{eq 1:proof Pro 1}
        \ell(x_0,h_{\gamma}^{i+1}(x_0))+\gamma W^{i}_{\gamma}(\hat{x}^{i+1}_{1})\leq \hat{V}^{i}_{\gamma}({x}_{0}),
    \end{equation}
    where $\hat{x}^{i+1}_{k+1}:=\hat{\phi}(1,x_k,h_{\gamma}^{i+1})$, $k\in\mathbb{Z}_{\geq 0}$.
    Also, since $\gamma\neq 0$, evaluating \eqref{eq3 Assumption for IPI} at $\hat{x}_1=\hat{\phi}(1,x_0,h_{\gamma}^{i+1}(x_0))\in\mathbb{R}^{n_x}$ leads to
    \begin{equation}\label{eq 2:proof Pro 1}
        \gamma\ell(\hat{x}^{i+1}_1,h_{\gamma}^{i+1}(\hat{x}_1))+\gamma^2 W_{\gamma}^{i}(\hat{x}^{i+1}_{2})\leq \gamma \hat{V}^{i}(\hat{x}_1).
    \end{equation}
    Using \eqref{eq 2:proof Pro 1} in \eqref{eq 1:proof Pro 1} yields 
    \begin{equation*}
        \ell(x_0,h^{i+1}_{\gamma}(x_0))+\gamma\ell(\hat{x}_1,h_{\gamma}^{i+1}(\hat{x}_1))+\gamma^2 W^{i}_{\gamma}(\hat{x}^{i+1}_{2})\leq \hat{V}^{i}_{\gamma}({x}_{0}).
    \end{equation*}
    Repeating this process for $N-2$ more times leads to
    \begin{equation}\label{eq 3:proof Pro 1}
        \sum_{k=0}^{N-1}\gamma^k\ell(\hat{x}^{i+1}_k,h_{\gamma}^{i+1}(\hat{x}_k))+\gamma^{N} W^{i}_{\gamma}(\hat{x}^{i+1}_{N})\leq \hat{V}^i_{\gamma}(x_0),
    \end{equation}
    where $\hat{x}^{i+1}_0=x_0$.
    Letting $N\rightarrow\infty$ and given $\hat{V}^i_{\gamma}(x_0)\geq 0$, $\forall x_0\in\mathbb{R}^{n_x}$, which hence can be dropped from the left hand side, inequality \eqref{eq 3:proof Pro 1} leads to $\hat{V}_{\gamma}^{i+1}(x)\leq\hat{V}_{\gamma}^{i}(x), \forall x\in\Omega$.
\hfill$\blacksquare$

The next proposition shows that the optimal policy for the incremental model verifies a $\mathcal{KL}$-stability property with respect to $\sigma$.
\begin{proposition}\label{pro: KL stability}
 If there exists $0<\gamma^{\star}\leq1$ such that $(1-\gamma^{\star})\bar{\alpha}_{V}(s)\leq\alpha_\Gamma(s), \forall s\in\mathbb{R}_{>0}$, then for any $\gamma\in(\gamma^{\star},\gamma_0)$, system~\eqref{eq:closed deffirentail inclusion} with optimal policies $h^{\star}_{\gamma}$ is $\mathcal{KL}$ stable with respect to $\sigma$,  i.e.,  there exists $\beta^{\star}\in\mathcal{KL}$ such that for any $x\in\Omega$, any solution $\hat{\phi}^{\star}(\cdot,x)$ to \eqref{eq:closed deffirentail inclusion} satisfies
 \begin{equation*}
  \sigma(\hat{\phi}^{\star}(k,x))\leq \beta^{\star}(\sigma(x),k), \forall k\in\mathbb{Z}_{\geq 0}.   
 \end{equation*}
 \hfill$\square$
\end{proposition}

\textit{Proof.}
    Let $\gamma\in(\gamma^{\star},\gamma_0)$, $x\in\Omega$, and $v=\hat{\phi}(1,x,h^{\star}_{\gamma}(x))$ with $h^{\star}_{\gamma}(x)\in H^{\star}_{\gamma}(x)$. 
    Since $\ell$ is non-negative and by using \eqref{eq1 Assumption for IPI}, it follows that
    \begin{equation*}
        \ell(x,h^{\star}_{\gamma}(x))\leq \hat{V}^{\star}_{\gamma}(x)\leq \hat{V}^{0}_{\gamma}(x)\leq \bar{\alpha}_{V}(\sigma(x)).
    \end{equation*}
    By definition of $\hat{V}^{\star}_{\gamma}$, one has that
    \begin{equation*}
        \hat{V}^{\star}_{\gamma}(x)=\ell(x,h^{\star}_{\gamma}(x))+\gamma \hat{V}^{\star}_{\gamma}(v),
    \end{equation*}
    therefore
    \begin{equation*}
        \hat{V}^{\star}_{\gamma}(v)-\hat{V}^{\star}_{\gamma}(x)=-\frac{1}{\gamma}\ell (x,h^{\star}_{\gamma}(x))+\frac{1-\gamma}{\gamma}\hat{V}^{\star}_{\gamma}(x),
    \end{equation*}
    which gives that
    \begin{equation}\label{eq1 in proof for optimal KL}
        \hat{V}^{\star}_{\gamma}(v)-\hat{V}^{\star}_{\gamma}(x)\leq-\frac{1}{\gamma}\ell(x,h^{\star}_{\gamma}(x))+\frac{1-\gamma}{\gamma}\bar{\alpha}_{V}(\sigma(x)).
    \end{equation}

    Define $\Upsilon_{\gamma}^{\star}:=\hat{V}^{\star}_{\gamma}+\frac{1}{\gamma}\Gamma$. 
    Combining \eqref{eq2 Assumption for IPI}  with \eqref{eq1 in proof for optimal KL} yields the following bounds
   \begin{equation*}
\alpha_\Gamma(\sigma(x))\leq\Upsilon^{\star}_{\gamma}(x)\leq \bar{\alpha}_{V}(x)+\frac{1}{\gamma^{\star}}\bar{\alpha}_W(\sigma(x)).
   \end{equation*}
   Denote $\bar{\alpha}_{\Upsilon}:=\bar{\alpha}_{V}+\frac{1}{\gamma^{\star}}\bar{\alpha}_W$ and $\underline{\alpha}_{\Upsilon}:=\alpha_\Gamma$.
   Moreover, from \eqref{eq1 in proof for optimal KL}, it follows that
   \begin{equation*}
       \Upsilon^{\star}_{\gamma}(v)-\Upsilon^{\star}_{\gamma}(x)\leq \frac{1}{\gamma}\big(-\alpha_\Gamma(\sigma(x))+(1-\gamma)\bar{\alpha}_{V}(\sigma(x))\big).
   \end{equation*}
   Since $(1-\gamma^{\star})\bar{\alpha}_{V}(s)\leq\alpha_\Gamma(s), \forall s\in\mathbb{R}_{>0}$ and $\gamma\geq \gamma^{\star}$, we have that
   \begin{equation*}
       \Upsilon^{\star}_{\gamma}(v)\leq \Upsilon^{\star}_{\gamma}(x)-\frac{1}{\gamma}\alpha_{\Upsilon}\big(\bar{\alpha}^{-1}_{\Upsilon}(\Upsilon^{\star}_{\gamma}(x)),\gamma\big),
   \end{equation*}
   where $\alpha_{\Upsilon}(\cdot,\gamma):=\frac{\gamma-\gamma^{\star}}{1-\gamma^{\star}}\alpha_\Gamma(\cdot)\in\mathcal{K}_{\infty}$.
   By induction, it can be seen that there exists $\beta^{\star}\in\mathcal{KL}$, such that
   \begin{equation*}
       \sigma(\hat{\phi}^{\star}_{\gamma}(k,x))\leq\beta^{\star}(\sigma(x),k), 
   \end{equation*}
   with $\beta^{\star}(s,k)\mapsto \underline{\alpha}_{\Upsilon}^{-1}(\{\min(\bar{\alpha}_{\Upsilon}(s),\gamma)\}^{k},k)$.
   This completes the proof.
\hfill$\blacksquare$

We are now give the results about the near optimality.
\begin{theorem}\label{th1}
  For any $x\in\Omega$, $i\in\mathbb{Z}_{\geq 0}$, $\gamma\in(\gamma^{\star},\gamma_0)$, and any solution to system~\eqref{eq:closed loop system}, 
  \begin{equation}\label{eq:condition in the1}
      \begin{aligned}
          \hat{V}^i_{\gamma}(x)-V^{\star}_{\gamma}(x)
          &\leq \bar{\alpha}_V(\beta^{\star}(\sigma(x),i),\gamma) +\frac{\gamma L_{\ell}\varepsilon_{\mathrm{IME}}}{1-\gamma}
      \end{aligned}
  \end{equation}
  with $\beta^{\star}\in\mathcal{K}\mathcal{L}$ form Proposition~\ref{pro: KL stability} and $\bar{\alpha}_V$ from Assumption~\ref{as:IPI}.
  \hfill$\square$
\end{theorem}

\textit{Proof.}
  Let $x\in\Omega$, $i\in\mathbb{Z}_{>0}$, $h^{i}_{\gamma}\in H^i_{\gamma}, h^{\star}_{\gamma}\in H^{\star}_{\gamma}$ and $\gamma\in(\gamma^{\star},\gamma_0)$.
   By Bellman equation and the definition of $\hat{V}^i_{\gamma}(x)$,
   \begin{equation}\label{eq1: proof of the1}
   \begin{aligned}
       \hat{V}_{\gamma}^i(x)-V_{\gamma}^{\star}(x)&=\hat{V}_{\gamma}^i(x)-\hat{V}_{\gamma}^{\star}(x)+\hat{V}_{\gamma}^{\star}(x)-V_{\gamma}^{\star}(x)\\
       & \leq \hat{V}_{\gamma}^i(x)-\hat{V}_{\gamma}^{\star}(x) +\frac{\gamma L_{\ell}\varepsilon_{\mathrm{IME}}}{1-\gamma},
   \end{aligned}
   \end{equation}
   which is deduced by using \eqref{eq:bound of incremental iotimal}.
   Since 
    \begin{equation*}
   h_{\gamma}^{i}(x)\in H_{\gamma}^{i}(x)=\arg\min_{u\in\mathcal{U}}\big\{\ell(x,u)+\gamma \hat{V}^{i-1}_{\gamma}(\hat{\phi}(1,x,u))\big\},
   \end{equation*}
   it follows that, $\forall u\in\mathcal{U}$,
   \begin{equation*}
        \begin{aligned}
       \hat{V}_{\gamma}^i(x)-\hat{V}_{\gamma}^{\star}(x)\leq& \ell(x,u)+\gamma \hat{V}^{i}_{\gamma}(\hat{\phi}(1,x,u))\\
       &-\ell(x,h^{\star}_{\gamma}(x))-\gamma \hat{V}^{\star}_{\gamma}(\hat{\phi}(1,x,h^{\star}_{\gamma}(x))).
   \end{aligned}
   \end{equation*}
   Therefore, by taking $u=\hat{h}^{\star}_{\gamma}\in \hat{H}^{\star}_{\gamma}(x)$ it follows that
   \begin{equation*}\label{eq2: proof of the1}
   \begin{aligned}
       \hat{V}^i_{\gamma}(x)-\hat{V}_{\gamma}^{\star}(x) &\leq \gamma\big(\hat{V}^{i-1}_{\gamma}-\hat{V}^{\star}_{\gamma}\big)(\hat{\phi}(1,x,h^{\star}_{\gamma})).
    \end{aligned}
   \end{equation*}
  Using Proposition~\ref{pro:improvement property} and repeating the above reasoning $i-1$ times, we obtain
  \begin{equation}\label{eq3: proof of the1}
   \begin{aligned}
       \hat{V}^i_{\gamma}(x)-\hat{V}_{\gamma}^{\star}(x) &\leq \gamma^i\big(\hat{V}^{0}_{\gamma}-\hat{V}^{\star}_{\gamma}\big)(\hat{\phi}(i,x,h^{\star}_{\gamma})).
    \end{aligned}
   \end{equation}
   Since $V^{\star}_{\gamma}\geq 0$, by Assumption~\ref{as:IPI}, we have
   \begin{equation}\label{eq4: proof of the1}
       \hat{V}^0_{\gamma}(x)-\hat{V}^{\star}_{\gamma}(x)\leq \hat{V}^0_{\gamma}(x)\leq \bar{\alpha}_V(\sigma(x),\gamma).
   \end{equation}
   Combining \eqref{eq3: proof of the1}, \eqref{eq4: proof of the1} with \eqref{eq1: proof of the1}, Proposition~\ref{pro: KL stability}, and noticing that $\bar{\alpha}_V$ is non-decreasing, we finally have \eqref{eq:condition in the1}.
\hfill$\blacksquare$

\subsection{Robust stability}
Before giving the robust stability results, we establish the following Lyapunov property for the system during the policy iteration process.
\begin{proposition}\label{pro:Lyapunov incremental}
  There exist $\underline{\alpha}_{Y}\in\mathcal{K}_{\infty}, \bar{\alpha}_Y,\alpha_Y:\mathbb{R}_{\geq0}\times(\gamma^{\star},\gamma_0)\rightarrow\mathbb{R}_{\geq0}$ of class $\mathcal{K}_{\infty}$ in their first argument such that for any $i\in\mathbb{Z}_{\geq0}$ there exist $Y^{i}_{\gamma}:\mathbb{R}^{n_x}\rightarrow\mathbb{R}_{\geq0}$ satisfying 
  \begin{itemize}
      \item[(i)] For any $x_k\in\mathbb{R}^{n_x}$, $\underline{\alpha}_{Y}(\sigma(x_k))\leq Y^i_{\gamma}(x_k)\leq \bar{\alpha}_{Y}(\sigma(x_k),\gamma)$;
      \item[(ii)] For any $x_k\in\mathbb{R}^{n_x}$, $Y^i_{\gamma}(\hat{x}_{k+1})-Y^i_{\gamma}({x}_{k})\leq \frac{1}{\gamma}\big(-\alpha_Y(\sigma(x_k),\gamma)+\Upsilon^i(\sigma(x_k),\gamma)\big)$,
  \end{itemize}
  for any $\gamma\in(\gamma^{\star},\gamma_0)$, where $\Upsilon^{i}:\mathbb{R}_{\geq0}\times(\gamma^{\star},\gamma_0)\rightarrow\mathbb{R}_{\geq0}$ is of class $\mathcal{K}_{\infty}$ in its first argument defined as $\Upsilon^{i}(\sigma,\gamma):=(1-\gamma)\gamma^i\bar{\alpha}_{V}(\beta^{\star}_{\gamma}(\sigma,i))$.
  \hfill$\square$
\end{proposition}

\textit{Proof.}
    Since the stage cost $\ell$ is non-negative, for all $x_k\in\mathbb{R}^{n_x}$, one has that
    \begin{equation*}
        \ell(x_k,h^i_{\gamma}(x_k))\leq \hat{V}^i_{\gamma}(x_k).
    \end{equation*}
Moreover, according to (PI.3)-(PI.6) in Algorithm~\ref{al PI-IADP}, $\hat{V}_{\gamma}^{i}(x_k)= W_{\gamma}^{i}(x_k)$, which yields that
\begin{equation}\label{eq1 proof pro2}
\begin{aligned}
    W^{i}_{\gamma}(\hat{x}_{k+1})&-W_{\gamma}^{i}(x_k)= W^{i}_{\gamma}(\hat{x}_{k+1})-\hat{V}^{i}_{\gamma}(x_k)\\
    =&W^{i}_{\gamma}(\hat{x}_{k+1})-\gamma W^i_{\gamma}(\hat{x}_{i+1})-\ell(x_k,u^i(x_k))\\
    =&(1-\gamma)W^{i}_{\gamma}(\hat{x}_{k+1})-\ell(x_k,u^i(x_k)).
\end{aligned}
\end{equation}
    Defining $Y^{i}_{\gamma}:=\hat{V}^i_{\gamma}+\frac{1}{\gamma}\Gamma$, it follows from \eqref{eq2 Assumption for IPI} that 
    \begin{equation*}
    \begin{aligned}
        Y^{i}_{\gamma}(x_k)\geq& \ell(x_k,h^i_{\gamma}(x_k))+\frac{1}{\gamma}\Gamma(x_k)\\
        &\geq \alpha_\Gamma(\sigma(x_k))=:\underline{\alpha}_{Y}(\sigma(x_k)).
    \end{aligned}
    \end{equation*}
    By using Assumptions~\ref{as:IPI} and
    Proposition~\ref{pro:improvement property}, one has 
    \begin{equation*}
    \begin{aligned}
        Y^i_{\gamma}(x_k)&\leq  \bar{\alpha}_V(\sigma(x_k),\gamma)\!+\!\frac{1}{\gamma}\bar{\alpha}_W(x_k)=:\bar{\alpha}_Y(\sigma(x_k),\gamma).        
    \end{aligned}
    \end{equation*}
    Therefore, item (i) is proven.

    From (PI.3)-(PI.6), it follows that
    \begin{equation}\label{eq1: pro 2}
        \begin{aligned}
            \hat{V}_{\gamma}^i(\hat{x}_{k+1})-\hat{V}_{\gamma}^i(x_k)=-\frac{1}{\gamma}\ell(x_k,h^i_{\gamma}(x_k))+\frac{1-\gamma}{\gamma}\hat{V}_{\gamma}^i(\hat{x}_{k+1})
        \end{aligned}
    \end{equation}
    In view of \eqref{eq2 Assumption for IPI} 
    and \eqref{eq1: pro 2}, we have
    \begin{equation*}
        \begin{aligned}
            Y_{\gamma}^i(\hat{x}_{k+1})-Y_{\gamma}^i(x_k)&\leq \frac{1-\gamma}{\gamma}(\hat{V}_{\gamma}^i(x_k)-\hat{V}_{\gamma}^{\star}(x_k))\\&-\frac{1}{\gamma}\alpha_\Gamma(\sigma(x_k))
            +\frac{1-\gamma}{\gamma}\hat{V}^{\star}_{\gamma}(x_k).
        \end{aligned}
    \end{equation*}
    By using Theorem~\ref{th1} and Proposition~\ref{pro: KL stability}, 
    \begin{equation}\label{eq3: proof Pro2}
    \begin{aligned}
        Y^i_{\gamma}(\hat{x}_{k+1})-Y^i_{\gamma}&(x_k)\leq\frac{1-\gamma}{\gamma}\gamma^i\bar{\alpha}_{V}(\beta_{\gamma}^{\star}(\sigma(x_k),i),\gamma)\\
        &-\alpha_\Gamma(\sigma(x))+(1-\gamma)\bar{\alpha}_{V^{\star}}(\sigma(x_k)).
    \end{aligned}
    \end{equation}
Since $(1-\gamma^{\star})\bar{\alpha}_{V^{\star}}(s)\leq \alpha_\Gamma(s), \forall s\in\mathbb{R}_{>0}$ as stated in Proposition~\ref{pro: KL stability}, the last two terms in \eqref{eq3: proof Pro2} satisfy 
\begin{equation*}
\begin{aligned}
   -\alpha_\Gamma(\sigma(x))+(1-\gamma)\bar{\alpha}_{V^{\star}}(\sigma(x_k))&\leq \frac{\gamma-\gamma^{\star}}{1-\gamma^{\star}}\alpha_\Gamma(x_k)\\
   =:\alpha_Y(x_k,\gamma). 
\end{aligned}
\end{equation*}
Define $\Upsilon^{i}(\sigma,\gamma):=(1-\gamma)\gamma^i\bar{\alpha}_{V}(\beta^{\star}_{\gamma}(\sigma,i))$. 
Item (ii) can be deduced and this completes the proof.
\hfill$\blacksquare$

According to the form of $\Upsilon_{\gamma}^i$, it follows that item~(ii) in  Proposition~\ref{pro:Lyapunov incremental} is a dissipative inequality of system~\eqref{eq:closed deffirentail inclusion} for which the supply rate consists of a negative term, namely $-\frac{1}{\gamma}\alpha_{Y}(\cdot,\gamma)$, and a non-negative term $\frac{1}{\gamma}\Upsilon^{i}(\cdot,\gamma)$ that can be made as small as desired by increasing $i$. 
Then, the following robust stability result is derived.

\begin{theorem}\label{th2}
Use the Lyapunov function definition and notation from Proposition~\ref{pro:Lyapunov incremental}.
For any $x\in\Omega$, given $\delta\geq 0$ and $\tilde{\delta}:=\underline{\alpha}_Y(\delta)>0$,   when Assumptions~\ref{as:existense of optimal control}-\ref{as:IPI} hold, there exists  $i^{\star}\in\mathbb{Z}_{\geq0}$, 
such that
\begin{equation}\label{eq:condition in th2}
    i^{\star}\geq \frac{\ln{\big(\frac{\alpha_Y(\bar{\alpha}_Y^{-1}(\underline{\alpha}(\delta),\gamma),\gamma)}{2(1-\gamma)\bar{\alpha}_V(\beta^{\star}(\underline{\alpha}_Y^{-1}(\bar{\alpha}_Y(\Delta,\gamma)),0),\gamma)}\big)}}{\ln{(\gamma)}},
\end{equation}
for any $i\geq i^{\star}$, system~\eqref{eq:closed loop system} is robustly $\mathcal{KL}$-stable. 
\hfill$\square$
\end{theorem}

\textit{Proof.}
Denote $\Delta:=\max\limits_{\rho\in\varepsilon_{\mathrm{IME}}\mathbb{B}}\big(\sigma(x+\rho)-\sigma(x)\big)\geq \sigma(\hat{\phi}(1,x,h^i_{\gamma}(x)))-\sigma(\phi(1,x,h^{i}_{\gamma}(x)))$,
define $\tilde{\Delta}_{\gamma}:=\bar{\alpha}_Y(\Delta,\gamma)>0$.
Using item (ii) in Proposition~\ref{pro:Lyapunov incremental}, define $v=\hat{\phi}(1,x,h^i_{\gamma}(x))$, one has that
\begin{equation}\label{eq1 in th2}
    Y_{\gamma}^i(v)-Y_{\gamma}^i(x)\leq \frac{1}{\gamma}\big(-\alpha_Y(\sigma(x),\gamma)+\Upsilon^i(\sigma(x),\gamma)\big).
\end{equation}
As $\Upsilon^i(\cdot,\gamma)$ is non-decreasing and $\alpha_Y(\cdot,\gamma)\in\mathcal{K}_{\infty}$, using item (i) of Proposition~\ref{pro:Lyapunov incremental} and the fact that $Y^i_{\gamma}(x)\leq \tilde{\Delta}_{\gamma}$, \eqref{eq1 in th2} yields
\begin{equation*}
    \begin{aligned}
        Y_{\gamma}^i(v)-Y_{\gamma}^i(x)\leq&\frac{1}{\gamma}\big(\Upsilon^i(\underline{\alpha}_Y^{-1}(\tilde{\Delta}_{\gamma})),\gamma\big)\\
        &-\alpha_Y(\bar{\alpha}_Y^{-1}(Y_{\gamma}^i(x),\gamma),\gamma).
    \end{aligned}
\end{equation*}
As $\beta^{\star}\in\mathcal{KL}$, for any $s\in\mathbb{R}_{\geq0}$ and $i\in\mathbb{Z}_{\geq 0}$, it follows that
\begin{equation*}
    \beta^{\star}(s,i)\leq\beta^{\star}(s,0).
\end{equation*}
As a result, when selecting $i^{\star}$ satisfying~\eqref{eq:condition in th2},
it follows that 
\begin{equation*}
\begin{aligned}
    \Upsilon^{i^{\star}}(\bar{\alpha}_Y^{-1}(\tilde{\Delta}_{\gamma}),\gamma)&\leq (1-\gamma)\gamma^{i^{\star}}\bar{\alpha}_V(\beta^{\star}(\underline{\alpha}_Y^{-1}(\tilde{\Delta}_{\gamma}),0),\gamma)\\
    &\leq\frac{1}{2}\alpha_Y(\bar{\alpha}_Y^{-1}(\tilde{\delta},\gamma),\gamma).
\end{aligned}
\end{equation*}
Consequently, for any $i\geq i^{\star}$, when $Y_{\gamma}^i\geq \tilde{\delta}$, 
\begin{equation*}
    Y_{\gamma}^i(v)-Y_{\gamma}^i(x)\leq-\frac{1}{2\gamma}\alpha_Y(\bar{\alpha}_Y^{-1}(Y_{\gamma}^i(x),\gamma),\gamma).
\end{equation*}
Therefore, there exists $\tilde{\beta}\in\mathcal{KL}$ such that for any solution $\hat{\phi}$ with respect to the incremental approximation model initialized at arbitrary $x$ and any $k\in\mathbb{Z}_{\geq 0}$,
\begin{equation*}
    Y^i_{\gamma}(\hat{\phi}(k,x,h^i_{\gamma}))\leq\max\{\tilde{\beta}(Y^i_{\gamma}(x),k),\tilde{\delta}\}.
\end{equation*}
Therefore, by using item (i) in Proposition~\ref{pro:Lyapunov incremental} and the definition of $\tilde{\delta}$, it follows that
\begin{equation}\label{eq 47 in proof th2}
    \sigma(\hat{\phi}(k,x,h^i_{\gamma}))\leq\max\{\beta(\sigma(x),k),\delta\},
\end{equation}
where $\beta(s,k):=\underline{\alpha}_Y^{-1}(\tilde{\beta}(\bar{\alpha}_Y(s,\gamma),k))$.
This concludes the proof.
\hfill$\blacksquare$

\begin{remark}
 Theorem~\ref{th2} shows that the proposed IPI framework guarantees robust $\mathcal{KL}$-stability, even without an initially stabilizing policy. This result is significant as it ensures that iterative policy updates naturally lead to stability despite model uncertainties and approximation errors. It also highlights the trade-off between convergence and robustness, as a smaller discount factor \(\gamma\) accelerates stability but may slow policy improvement, while a larger \(\gamma\) speeds up optimality convergence but requires more iterations for stability. This theorem strengthens IPI’s applicability in model-free control by ensuring stability in unknown nonlinear systems, making it a reliable approach for adaptive optimal control.  
 \hfill$\square$
\end{remark}

\section{Simulations}\label{sec:simulation}
Consider a nominal nonlinear system (Model A) govern by the followig dynamics
\begin{equation}\label{eq:simulation}
    x_{k+1}=\begin{bmatrix}
        x_{2,k}\\
        -2x_{1,k}-3x_{2,k}+\sin(x_{1,k})+u_k
    \end{bmatrix}
\end{equation}
where $x_{k}=[x_{1,k};x_{2,k}]\in\mathbb{R}^2, k\in\mathbb{Z}_{\geq 0}$.
Select positive definite matrix $Q\in\mathbb{R}^{2\times 2}$ and positive constant $R>0$.
The objective is to use the proposed IPI to steer the system from a given initial state \(x_0\) to the equilibrium $(0,0)$, while minimizing the following performance index
\begin{equation}\label{eq:cost in simulation}
    J=\sum_{k=0}^{\infty}\gamma^k \big(x_k^{\top} Q x_k+R u_k^2\big)
\end{equation}
under the condition of knowing only the system input \(u\) and state \(x\) at the current and previous time steps $t_k$ and $t_{k-1}$.

Here, we consider the value approximator $W^{i}_{\gamma}(x_k)$ as $x_k^{\top}P^{(i)}_{\gamma}x_k$, where $P^{(i)}_{\gamma}\in\mathbb{R}^{2\times 2}$ is a positive definite matrix to be found recursively.
In this case, from (PI.1)-(PI.6), one has that
\begin{equation}\label{eq:incremental simulation}
\begin{aligned}
    &\Delta u^i_k =-(R+\gamma \hat{B}_{k-1}^{\top} P^{(i)}_{\gamma}\hat{B}_{k-1})^{-1}\times\\
    &~~\big[Ru_{k-1}+\gamma\hat{B}_{k-1}^{\top} P^{(i)}_{\gamma} x_k+\gamma\hat{B}_{k-1}P^{(i)}_{\gamma}\hat{A}_{k-1}\Delta x_k\big].
\end{aligned}
\end{equation}
Therefore,  we can conclude that the policy is in the feedback form of
system variables $(u_{k-1},x_k,\Delta x_k)$, and the gains are function of the current incremental model $(\hat{A}_{k-1},\hat{B}_{k-1})$.

For simplicity, select $Q=I_2$ and $R=1$.
The implement procedure is illustrated as follows.

\begin{figure}[htb]
    \centering
    \includegraphics[width=0.9\linewidth]{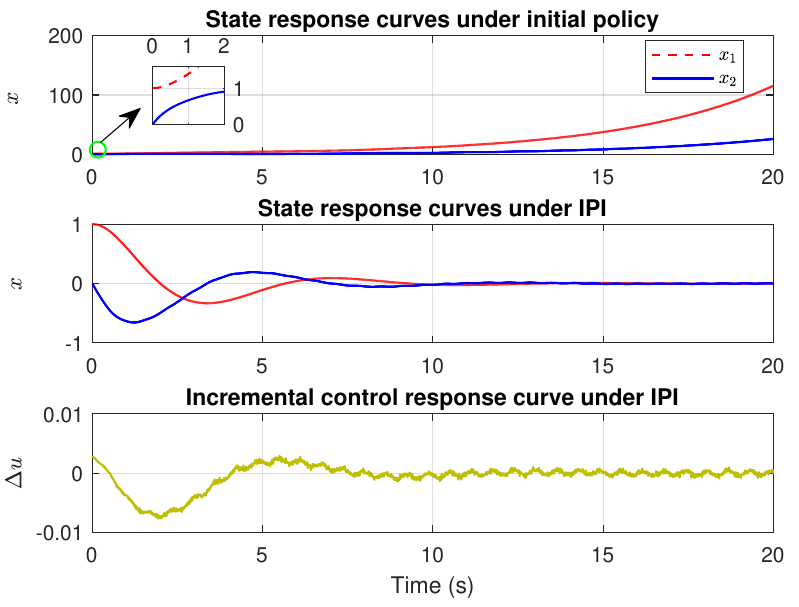}
    \caption{State and incremental control response curves of the proposed method.}
    \label{fig:response}
\end{figure}

\textit{Offline training:} Using \eqref{eq:simulation} (Model A), a set of data \((\mathbf{x}_{0:N}, \mathbf{u}_{0:N})\) is collected by applying a randomly generated input signal (a sinusoidal signal is considered in this paper). 
The corresponding \(\hat{A}_{k-1}, \hat{B}_{k-1}\) are identified by using batch LS for offline policy iteration.
To verify Theorems~\ref{th1}-\ref{th2}, during the offline training, choose an initial policy satisfying \eqref{eq1 Assumption for IPI}, $u^0_k=[-2.5\ -1]x_k$ and $\gamma=0.7$.
It can be seen from the upper subfigure in Fig.~\ref{fig:response} that the initial policy $u^{0}_{\gamma}$ is an unstable policy.
 $P^i_{\gamma}$ is updated by solving
\begin{equation}
    x_k^{\top}P_{\gamma}^{(i+1)}x_k = x_k^{\top}Qx_{k}+R (u_k^{i})^2+\gamma \hat{x}_{k+1}^{\top}P^{(i)}_{\gamma}\hat{x}_{k+1},
\end{equation}
which implies that there exists $\alpha_\Gamma$ satisfying $\eqref{eq2 Assumption for IPI}$.
Moreover, since $\Delta u^i_k$ is the analytic solution of (PI.5), it always makes \eqref{eq3 Assumption for IPI} hold.

\textit{Online implementation:} To verify the robustness of the proposed method, we consider a different physical system 
(Model B) to control:
\begin{equation}\label{eq:simulation 2}
    x_{k+1}\!=\!\begin{bmatrix}
        x_{2,k}\\
        -2x_{1,k}\!-\!0.5x_{2,k}\!+\!\sin(x_{1,k})\!+\!0.2u_k\!+\!u_{d,k}
    \end{bmatrix},
\end{equation}
where $u_{d,k}$ is the disturbance in the form of 
\begin{equation*}
    u_{d,k}=0.2 \sin(0.1 t_k) + 0.1 w(k)
\end{equation*}
with Gaussian noise $w(k)$.
The online iteration is implemented in a recursive manner with the offline trained policy as a baseline policy.
That is, the kernel matrix $P$ is updated for each time step $t_k$: 
\begin{equation}
    x_{k}^{\top}P_k x_k=x_{k}Qx_k+R u_{k}^2+\gamma\hat{x}_{k+1}^{\top}P_{k-1}\hat{x}_{k+1}.
\end{equation}
The incremental policy $\Delta u_k$ is improved based on the new $P_k$ by using \eqref{eq:incremental simulation}.
The state and incremental control response curves are shown in the second and third subfigures in Fig.~\ref{fig:response}.
Thus, the trained policy makes the system stable while rejecting the uncertainty and disturbances brought by Model~B.


\textit{Comparison:} To further demonstrate the robustness of the proposed IPI approach, we compare it with the traditional nonlinear PI-ADP method~\cite{Liu2013Policy}. 
In the conventional approach, a model neural network is first trained to approximate the system dynamics, and then this learned model is used to train the actor and critic networks via policy iteration with an initial stable policy. 
However, this indirect learning process introduces modeling errors that propagate through the control design, leading to degraded performance under uncertainties and disturbances when implementing online.

\begin{figure}[htb]
    \centering
    \includegraphics[width=0.9\linewidth]{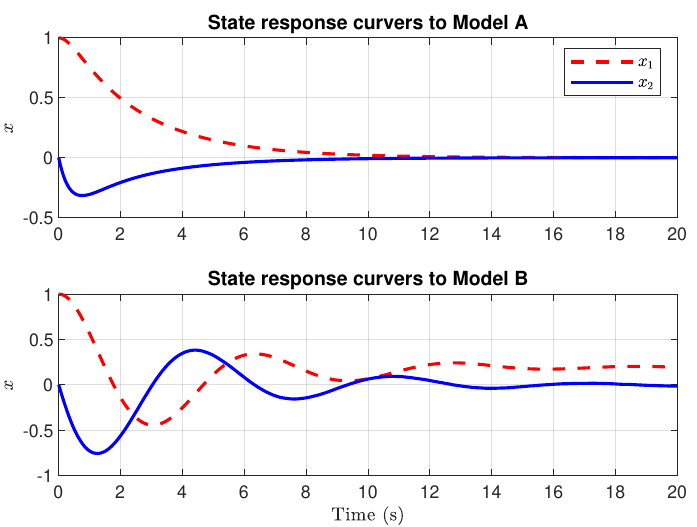}
    \caption{State response curves of the traditional PI-ADP method~\cite{Liu2013Policy}.}
    \label{fig:comparison}
\end{figure}

As both methods are trained using data collected from Model A, their ability to generalize to Model B under uncertainties and disturbances is evaluated.
Fig.~\ref{fig:response} demonstrates that the proposed IPI approach remains robust when applied to Model B, effectively handling disturbances and model uncertainties. In contrast, Fig.~\ref{fig:comparison} shows that the traditional PI-ADP method, despite being trained on Model A, fails to maintain stability when uncertainties and disturbances are introduced. 
This degradation highlights the limitations of relying on a pre-trained model for control design, since inaccuracies in the learned dynamics can negatively impact policy performance. 
The IPI approach mitigates these issues by continuously adapting its policy, ensuring superior robustness and stability across different operating conditions.

\section{Conclusions}\label{sec:conclusion}
In this paper, a general model-free incremental policy iteration framework for nonlinear systems is proposed, employing the recursive least squares method to identify linear approximation system matrices.
This allows the offline pre-trained policy to be updated online with limited data.
The approach avoids the high training cost and poor interpretability associated with the global approximation used in traditional nonlinear ADP methods, while robustly adapting to dynamic system variations through an incremental update mechanism.
The near-optimality and robust stability of the algorithm are theoretically proven, providing a solid theoretical foundation. 
Future work will focus on extending the method to continuous systems and its applications in engineering practice.

\vspace{12pt}

\end{document}